\documentstyle[12pt,amsfonts,amssymb,leqno]{article}

\newfont\got{eufm10}

\newcounter{secnum}

\begin{document}
\setcounter{section}{+1}

\begin{center}
{\Large \bf A note on an alleged proof of} 
\end{center}

\begin{center}
{\Large \bf 
the relative consistency of $P=NP$ with $PA$}
\end{center}

\begin{center}
\renewcommand{\thefootnote}{\fnsymbol{footnote}}
\renewcommand{\thefootnote}{arabic{footnote}}
\renewcommand{\thefootnote}{\fnsymbol{footnote}}
{\large Ralf-Dieter Schindler}
\renewcommand{\thefootnote}{arabic{footnote}}
\end{center}
\begin{center} 
{\footnotesize
{\it Institut f\"ur formale Logik, Universit\"at Wien, 1090 Wien, Austria}} 
\end{center}

\begin{center}
{\tt rds@logic.univie.ac.at}

{\tt http://www.logic.univie.ac.at/${}^\sim$rds/}\\
\end{center}


\bigskip
N.C.A. da Costa and F.A. Doria claim to have shown in \cite{cd} that $P=NP$ is
relatively consistent with $PA$. The purpose of the present note is to argue that
there is a mistake in that paper. Specifically, we want to point out that Corollary
5.14 of \cite{cd} -- which is used in the proof of their main result --
is probably false.

We are first going to reconstruct their argument. According to that reconstruction,
the argument of \cite{cd} would in fact show that $PA$ proves $P=NP$. We'll then
discuss Cor.
5.14 of \cite{cd}. 
However, rather than talking about provability in $PA$ or stronger theories, we'll
stick to a different attitude and argue internally\footnote{Argumenting externally 
might be more in line with the paper itself. However,
such a reconstruction typically starts from a listing of Turing machines which is
such that $PA$ won't be able to prove ``$P<NP \Leftrightarrow f$ is total.'' 
(The definition of $f$ is given in the main text.)}: 
using Cor.
5.14 of \cite{cd} we'll derive a contradiction from the assumption that $P<NP$ (and we
implicitly assume our argument goes thru in $PA$).

We'll follow the notation of \cite{cd} (with the exception of $f_{\lnot A}$).

\bigskip
{\em Some Turing machines.} $\bigvee(z) = 1$ iff $\pi_1(z)$ codes a cnf-Boolean
expression and $\pi_2(z)$ codes an assignment which satisfies it (o.w. $\bigvee(z) =
0$). $E$ is a fixed exponential Turing machine that solves any instance of the
satisfiability problem (in particular, $\bigvee(\langle z,E(z) \rangle) = 
1$ for any $z$ coding a
satisfiable cnf-Boolean
expression). Let $G$ be a partial recursive function; 
for $n<\omega$ such that $G(n) \downarrow$ we let $Q^{G(n)}$ 
be that Turing machine which, given an input $z$,
first computes 
$G(n)$ and then computes 
$E(z)$ in case $z \leq G(n)$ and $0$ in case $z > G(n)$
(in particular, 
$Q^{G(n)}(z) = E(z)$ for $z \leq G(n)$ and $Q^{G(n)}(z) = 0$ for $z > G(n)$; 
cf. p. 10 of \cite{cd}). Notice $\bigvee$ and all $Q^{G(n)}$ (for $n<\omega$
and $G(n) \downarrow$) are
polynomial time Turing machines (in fact the $Q^{G(n)}$'s are ``finite''). 
By the Baker-Gill-Solovay trick there is a recursive
enumeration $( P_m \colon m<\omega )$ of all polynomial Turing machines;
we may think of (the ``G\"odel number'' of) $P_m$ as an ordered pair $\langle
g , c \rangle$ where $g$ codes $P_m$'s program and $c$ is a code for a polynomial
clock.

{\em Some functions.} We let $f(m)$ be the least $z$ such that 
$\bigvee(z) = 1$, whereas $\bigvee(\langle \pi_1(z),P_m(\pi_1(z)) 
\rangle) = 0$ (i.e., $f(m)$
witnesses that $P_m$ doesn't prove $P=NP$). We have: $f$ is recursive, and
$f$ is total iff $P<NP$.
($f$ is written $f_{\lnot A}$ in \cite{cd}; cf. \cite{cd} p. 4.)
Let $(\psi_i \colon i<\omega)$ be a recursive enumeration of all polynomial functions 
from $\omega$ to $\omega$.
We let $F(m) = max \{ f \circ \psi_i(m) \colon m \leq i \wedge f \circ \psi_i(m)
\downarrow \} +1$. (We understand that $max \ \emptyset = 0$.) Note that $F$
dominates $f \circ \psi$ 
(in the sense that $F(m) > f \circ \psi(m)$ for all sufficiently large $m$ with
$f \circ \psi_i(m)
\downarrow$) for any polynomial $\psi$. If $f$ is total then $F$ is recursive. 

Corollary
5.14 of \cite{cd} now reads as follows: {\bf Main Lemma.}
If $F$ is recursive, then there is a {\em linear} $\psi \colon \omega
\rightarrow \omega$ s.t. for all $m$ and $n$ do we have that 
$Q^{F(m)}(n) = P_{\psi(m)}(n)$.
Let us also consider the following, which trivially follows from this Main Lemma:
{\bf Main Lemma'.} If $F$ is recursive, then 
there is a {\em polynomial} $\psi \colon \omega
\rightarrow \omega$ s.t. for all $m$ and $n$ do we have that 
$Q^{F(m)}(n) = P_{\psi(m)}(n)$.

Given the Main
Lemma' we may now prove $P=NP$ as follows. 
Suppose not. Then $f$ is total recursive, and hence so is $F$. If $\psi$ 
is as in the Main
Lemma'  
then $F(m) > f \circ \psi(m)$ for all sufficiently large $m$. On the other hand,
$f(\psi(m))$ is the least $z$ such that 
$\bigvee(z) = 1$, whereas $\bigvee(\langle 
\pi_1(z),Q^{F(m)}(\pi_1(z)) \rangle) = 0$. For $\pi_1(z)
\leq F(m)$ we'll have $Q^{F(m)}(\pi_1(z)) = E(\pi_1(z))$, so that $f(\psi(m)) \geq
\langle F(m)+1,0 \rangle \geq F(m)+1$. 
We'll thus have $F(m) > f
\circ \psi(m) \geq F(m)+1$ for all sufficiently large $m$. Contradiction! We have
shown that $P=NP$.

Have we? Not so, I claim. Let's discuss the Main Lemma. 
The Baker-Gill-
Solovay trick uses the device of ``clocks'' in order to arrive at a recurive
enumeration $(P_m \colon m<\omega)$ of all polynomial Turing machines. 
I.e. (cf. above), the G\"odel number of $Q^{F(m)}$, viewed as a machine with
a clock attached to it, will be -- typically -- the ordered pair $\langle g,c
\rangle = \langle g(m),c(m)
\rangle$ of a code $g$ for 
$Q^{F(m)}$'s program (without a clock) and a code $c$ for a polynomial clock.
Now the clock is not supposed to shut down 
$Q^{F(m)}$'s operation before $F(m)$ gets known. 
That is, $c$ depends on the length of the
computation of $F(m)$; in fact, if $F$ is ``complicated,'' $c \approx$
the length of the
computation of $F(m)$. Hence the function $m \mapsto$ G\"odel number $\langle 
g(m),c(m)
\rangle$ of
$Q^{F(m)}$ will be 
at least as ``complex'' as $m \mapsto$ the length of the computation of
$F(m)$.
There is no reason to believe that it should be linear (or, polynomial, for that matter).
In other words, any function 
$\psi$ s.t. 
$Q^{F(m)}(n) = P_{\psi(m)}(n)$ for all $m$ and $n$
will be as ``complex'' as $F$ is; which is,
I think, as it should be. 
But this then poses a serious problem. Suppose a $\psi$ as
above can only be as ``complex'' 
as $F$. The above argument breaks down without an $F$ s.t.
$F$ dominates $f \circ 
\psi$. 

It is hard to believe that there should be a proof of the Main
Lemma' which doesn't actually prove (without assuming $P<NP$) the following.
\boldmath $(*) \ $\unboldmath For any total recursive $G$ there is a
polynomial $\psi \colon \omega \rightarrow \omega$ s.t. for all $m$
and $n$ do we have that $Q^{G(m)}(n) =
P_{\psi(m)}(n)$. But \boldmath $(*) \ $\unboldmath is false, as we shall now show. 
Recall that $Q^{G(m)}$ always first computes $G(m)$; in this sense, the value $G(m)$
``shows up (on the tape) during the 
calculation of $Q^{G(m)}(n)$,''\footnote{We leave it
to the reader to make this precise.} for every $m$ and
$n$. Now let us consider $$G(m) = 
max \{ t \ \colon \ t {\rm \ shows \ up \ during \ the
\ calculation \ of \ } P_i(m) {\rm \ for \ } i \leq m^m \} + 1.$$
Clearly, $G$ is total recursive. Suppose that $\psi$ is as in 
\boldmath $(*) \ $\unboldmath for $G$. In particular, for every $m$,
$G(m)$ shows up during the 
calculation of $Q^{G(m)}(m)$. However, pick $m$ s.t. $\psi(m) \leq m^m$.
Then by construction, $G(m) > t$ for all $t$ which show up in the calculation of
$P_{\psi(m)}(m) = Q^{G(m)}(m)$. Contradiction!
\footnote{Another way to look at \cite{cd} is the following.
If the argument of \cite{cd} worked, we could re-run it by talking about finite
time machines rather than polynomial time machines (by finite time I
mean no matter how long the input is the machine will stop after c steps,
where c is a constant being independent from the length of the input). 
All $Q^{F(m)}$'s are finite time. So by starting 
from an enumeration of finite time machines rather than of polynomial
time machines the argument should really prove (the relative consistency
of) ``finite'' $= NP$ rather than just of $P = NP$, which truly 
is absurd. }  

I don't understand the alleged proof of Cor.
5.14 of \cite{cd} which appears on pp. 13 ff. of \cite{cd}.
Specifically, I don't understand Remark 5.11 on p. 14 of \cite{cd}. I think it
contains a statement (``Again we have, [...].'') which has not been verified.
The statement is reminiscent to 
\boldmath $(*) \ $\unboldmath above.

\bigskip
We believe that at some point s.o. will prove $P<NP$ by some form of ``Galois theory.''

\bigskip
I cordially thank N.C.A. da Costa and Chico Doria for their interest in my e-mail 
messages concerning \cite{cd}. 

\end{document}